\documentclass{iise}

\usepackage{multirow}
\usepackage{amsmath}
\usepackage{breqn}
\conference{Proceedings of the IISE Annual Conference \& Expo 2023 \\
K. Babski-Reeves, B. Eksioglu, D. Hampton, eds.}

\title{\titlesize Stochastic Optimization of Global Agrochemical Supply Chains with Risk Management: Modeling and Reformulation}

\author{
Saba Ghasemi Naraghi, Zheyu Jiang\\Oklahoma State University\\Stillwater, Oklahoma, USA 74078\\
\vspace{0.3cm}
}
\authorlist{Ghasemi Naraghi, Jiang} 
\abstractID{3611}

\DeclareMathOperator{\clean}{\text{clean}}
\DeclareMathOperator{\start}{\text{start}}
\DeclareMathOperator{\norm}{\text{norm}}

\DeclareMathOperator{\capacity}{\text{Cap}}
\DeclareMathOperator{\Inv}{\text{Inv}}
\DeclareMathOperator{\Demand}{\text{Demand}}
\DeclareMathOperator{\ship}{\text{ship}}
\DeclareMathOperator{\fixed}{\text{fixed}}
\DeclareMathOperator{\var}{\text{variable}}
\DeclareMathOperator{\hold}{\text{hold}}
\DeclareMathOperator{\loss}{\text{loss}}
\DeclareMathOperator{\Slack}{\text{Slack}}
\DeclareMathOperator{\variance}{\text{Variance}}

\begin{document}
\maketitle

\begin{abstract}
{\small The global agrochemical market is highly consolidated, with large multinational companies accounting for a major share of the market. Thus, even for a single agrochemical product, its global supply chain typically involves numerous paths connecting the raw material sources to the final customers. Besides structural complexity, agrochemical supply chains are also subject to seasonality and other unique uncertainties, thereby posing a need for risk management tools and strategies. In this study, we model and optimize an agrochemcial supply chain by developing and solving a stochastic mixed-integer quadratic constrained program (MIQCP). We model and control the demand uncertainty in this scenario-based MIQCP using variance. We also apply perspective reformulation techniques to convert the MIQCP to a mixed-integer linear program (MILP). Computational experiment results from an illustrative example show that, successively introducing perspective cuts to the reformulated MILP not only leads to a tight approximation of the original MIQCP model, but is also more computationally efficient than directly solving the MIQCP.}
\end{abstract}

\section*{Keywords}
Supply chain optimization, Agrochemical supply chain, Stochastic optimization, Perspective cuts, Mixed-integer quadratic constrained programming

\section{Introducttion}
In 2050, the global population is expected to increase to 9.7 billion, which puts unprecedented stress on food, energy, and water resources. Specifically, global food production must increase by at least 70\% between now and 2050 \cite{searchinger18}. To meet this growing demand in food, the production, processing, and distribution of agrochemicals, which include pesticides, herbicides, fungicides, and insecticides, need to be efficient and resilient. On the other hand, the global agrochemical market is highly consolidated, in which more than 60\% of global market share is owned the top four agrochemical conglomerates \cite{2017}. Each conglomerate owns a complex, diversified product line, whose associated supply chains are multi-stage networks involving numerous paths connecting the raw material sources eventually to the final customers. Besides structural complexity, agrochemical supply chains are further complicated by seasonality and various uncertainties caused by climate change, more frequent black swan events, and increasingly complex geopolitical landscapes across the world. To design efficient and resilient agrochemical supply chains, in this work, we develop an agrochemical supply chain optimization framework with risk management, and propose a reformulation strategy to efficiently solve optimization problem.

Among numerous recent works on supply chain optimization, Bassett and Gardner \cite{bassett09} presented a mixed-integer linear programming (MILP) formulations for global agrochemical supply chain optimization considering seasonality and uncertainties in customer demand. Liu and Papageorgiou \cite{lui12} extended the agrochemical supply chain optimization framework by modeling and comparing different plant expansion strategies. To further ensure continuous use and inactivity of warehouses for continuous periods of time, Brunaud et al. developed dynamic contract policy constraints for warehouses and incorporated them to the agrochemical supply chain model \cite{brunaud17}. In terms of incorporating uncertainty and risk quantification in supply chain optimization, You et al. \cite{you09} proposed a scenario-based two-stage stochastic linear programming framework and decomposition strategies for multi-product supply chain planning under demand and freight rate uncertainties. Later, Carneiro et al. \cite{carneiro} focused on the oil supply chain optimization problem, in which they incorporated Conditional Value-at-Risk (CVaR) as a risk assessment measure that quantifies the tail risk in their investment portfolio. Recently, we \cite{jiang} proposed a scenario-based two-stage mixed-integer nonlinear programming (MINLP) model for agrochemical supply chain optimization and adopt the concepts of Value-at-Risk (VaR) and CVaR to quantify and control the risks associated with demand unfulfillment. 

In this paper, we develop a scenario-based two-stage mixed-integer quadratic constrained programming (MIQCP) model for agrochemical supply chain optimization subject to demand uncertainties. Lim et al. \cite{lim11} discussed the limitations of using mean-CVaR measure in portfolio optimization, one of which is its sensitivity to data outliers, which can result in unstable and unpredictable outcomes. The mean-CVaR optimization problem is also nonconvex, which may impact its convergence to global optimality. Instead of using CVaR, in this work, we use variance to characterize the risks associated with demand unfulfillment. To linearize the nonlinear variance constraint, we introduce perspective cuts originally proposed by Frangioni and Gentile \cite{frangioni05} to reformulate the MIQCP model into a MILP. Frangioni and Gentile \cite{frangioni05} showed that the convex envelope of the objective function containing semicontinuous variables in a general mixed-integer program (MIP) is the perspective function of MIP's continuous part. Later, Bestuzheva et al. \cite{bestuzhev21} introduced the concept of perspective cuts to handle nonlinear constraints in MINLPs. Consider a MINLP with a linear objective function: $\min f(x,y,z)$, subject to $y \in \Omega,\, z \in \{0,1\}^n$, and nonlinear constraint $g(x) \leq 0$, in which $x \in \mathbb{R}^n$ are semi-continuous variables. In other words, for every $i=1,\cdots,n$, $x_i=0$ when $z_i=0$ and $x_i \in [l,u]$ when $z_i=1$. The perspective cuts to linearize the nonlinear constraint $g(x) \leq 0$ are given by: 

\begin{equation}
    \langle \nabla g(\bar{x}^{(k)}) , x^{(k)} \rangle + \left( g(\bar{x}^{(k)}) - \langle \nabla g(\bar{x}^{(k)}) , \bar{x}^{(k)}\rangle \right) z \leq 0, \quad \forall k \in K,
\end{equation}
where $\bar{x}^{(k)} = x^{k-1}$ is the solution vector of the problem in iteration $k-1$ (after adding $k-1$ cuts). Note that $\bar{x}^{(1)}_i$ is an arbitrary value in $[l_i,u_i]$. After replacing $g(x) \leq 0$ with these perspective cuts, the MINLP is reformulated to a MILP, which can be solved iteratively. Bestuzheva et al. (2021) also conducted a detailed computational study of perspective reformulation for MINLPs with convex and nonconvex nonlinear constraints. They showed that the perspective reformulation of convex MINLPs provides much tighter approximation of the original problems compared to conventional branch-and-cut approaches, thereby leading to significant computational time reduction \cite{bestuzhev21}. 

\section{Supply Chain Problem Statement}
We consider the supply chain of an agrochemical active ingredient (AI) produced in batch mode at one company. Each batch of AI takes $L$ time periods (e.g., days) to produce. AI production plant can only produce one batch at a time. The capacity of each batch is typically fixed by the capacities of the reactor and downstream processing units, such as distillation columns, crystallizers, and dryers. However, there is room for slight increase of batch capacity in order to satisfy the demand under some special occasions. After manufacturing $B$ batches of AI in each AI production plant, the plant must be inactive and gets cleaned for $K$ time periods. The produced AI is then sent to the warehouses/distribution centers, where it is distributed to different market regions. All facilities (AI production plants and warehouses) are connected by transportation links. Each market region is served by at least one distribution center connected via a transportation link. For each transportation link, there can be more than one mode of transportation (e.g., small vs. large truck, truck vs. rail). For instance, if 1.40 metric tons of AI product needs to be shipped, we can choose to use either two small trucks (each has a capacity of 0.7 metric ton) or one large truck (with a capacity of 1.5 tons). Each AI production plant knows the AI product demand for the first planning period (e.g., first week). However, starting from the second planning period, the AI product demand becomes uncertain. The distribution of the weekly demand can be modeled by fitting historical demand data. One of the requirements to ensure resilience of supply chain is to keep the weekly demand loss at a low level subject to different scenarios and uncertainties.

\section{Problem Formulation}
\subsection{MIQCP Formulation}
Following the problem statement, we present the following objective function for our MIQCP model. For every AI production plant $i \in I$, we are given the initial inventory level $\Inv_{i,t=0,s}$ for every demand scenario $s \in S$, all cost parameters (including inventory holding cost $c_i^{\hold}$, AI production fixed cost $c_i^{\fixed}$ and variable cost $c_i^{\var}$). Similarly, for every warehouse/distribution center $j \in J$, we are given its holding cost $d_j^{\hold}$ and fixed cost $d_j^{\fixed}$. For each transportation mode $m \in M$, we are given its shipping cost $c_{m}^{\ship}$. Finally, for each (weekly) planning period $t' \in T'$, the cost due to demand loss is given by $r_{t'}^{\loss}$. Hence, the objective function for this MIQCP is:

\begin{equation}
    \begin{split}
        \min & \sum_{s \in S} \rho_s (\sum_{t \in T} (\sum_{i \in I} (c_i^{\fixed} \alpha_{i,t,s} + c_i^{\var} x_{i,t,s} + c_i^{\hold} \Inv_{i,t,s} + \sum_{m \in M} \sum_{j \in J} c_m^{\ship} u_{i,j,m,t,s})\\
        & + \sum_{j \in J} (d_j^{\fixed} \beta_{j,t,s} + d_j^{\hold} w_{j,t,s} + \sum_{m \in M} \sum_{k \in K} c_m^{\ship} u_{j,k,m,t,s})) + \sum_{t' \in T'} r_{t'}^{\loss} \text{Loss}_{t',s})),
    \end{split}
\end{equation}
where $\alpha_{i,t,s}$ (resp. $\beta_{j,t,s}$) is a binary variable indicating whether plant $i$ (resp. warehouse $j$) is active in production (resp. active in storage) ($=1$) at time period $t$ under scenario $s$ or not ($=0$), continuous variable $x_{i,t,s}$ (resp. $w_{j,t,s}$) represents for the amount of AI that is being produced (resp. stored) in plant $i$ (resp. in warehouse $i$) at time period $t$ under scenario $s$, continuous variable $u_{i,j,m,t,s}$ denotes the amount of AI product that is shipped from $i$ to $j$ using transportation mode $m$ during time period $t$ under scenario $s$, and $\text{Loss}_{t',s}$ is the amount of demand loss for week $t'$ under scenario $s$. The constraints for our MIQCP are:

\begin{align}
    \text{s.t. } & \capacity^{\norm}_i \alpha_{i,t,s}^{\start} \leq x_{i,t+L-1,s} \leq \capacity^{\max}_i \alpha_{i,t,s}^{\start}, \quad \forall i \in I, t \in T, s \in S,\\
    & \alpha_{i,t,s}^{\start} - \alpha_{i,\tau,s} \leq 0, \quad \forall i \in I, t \in T, \tau \in \{t, \cdots, t+L-1\}, s \in S, \\
    & \alpha_{i,\tau,s}^{\start} \leq 1 - \alpha_{i,t,s}^{\start}, \quad \forall i \in I, t \in T, \tau \in \{t+1, \cdots, t+L -1\}, s \in S, \\
    & \frac{\sum_{\tau=1}^{t} \alpha^{\start}_{i,\tau,s}}{B_i} - 1 + \frac{1}{B_i}\leq \sum_{\tau=1}^{t+L} \alpha_{i,\tau,s}^{\text{clean}} \leq \frac{\sum_{\tau=1}^{t} \alpha^{\start}_{i,\tau,s}}{B_i}, \quad \forall i \in I, t \in T (t+L \leq |T|), s \in S, \\
    & \alpha_{i,t,s}^{\text{clean}} \leq 1 - \alpha^{\start}_{i,t,s}, \quad \forall i \in I, t \in T, s \in S,\\
    & \alpha_{i,\tau,s} \leq 1 - \alpha_{i,t,s}^{\text{clean}}, \quad \forall i \in I, t \in T, \tau \in \{t, \cdots, t+C-1\}, s \in S,\\
    &  \Inv_{i,t,s} + \sum_{j \in J} y_{i,j,t,s} = \Inv_{i,t-1,s} + x_{i,t,s}, \quad \forall i \in I, t \in T, s \in S, \\
    & \capacity^{\min,\ship}_{m} \lambda_{i,j,m,t,s} u_{i,j,m,t,s} \leq o_{i,j,m,t,s} \leq \capacity^{\max,\ship}_{m} \lambda_{i,j,m,t,s} u_{i,j,m,t,s}, \quad \forall i \in I, j \in J, m \in M, t \in T, s \in S, \\
    & \sum_{m \in M} o_{i,j,m,t,s} = y_{i,j,t,s}, \quad \forall i \in I, j \in J, t \in T, s \in S, \\
    & \sum_{m \in M} \lambda_{i,j,m,t,s} = 1, \quad \forall i \in I, j \in J, t \in T, s \in S, \\
    & 0 \leq w_{j,t,s} \leq \capacity^{\max}_j \beta_{j,t,s}, \quad \forall j \in J, t \in T, s \in S,\\
    & w_{j,t,s} + \sum_{k \in K} p_{j,k,t,s} = w_{j,t-1,s} + \sum_{i \in I} y_{i,j,t,s}, \quad \forall j \in J, t \in T, s \in S, \\
    & \capacity^{\min,\ship}_{m} \lambda_{j,k,m,t,s} u_{j,k,m,t,s} \leq o_{j,k,m,t,s} \leq \capacity^{\max,\ship}_{m} \lambda_{j,k,m,t,s} u_{j,k,m,t,s}, \quad \forall j \in J, k \in K, m \in M, t \in T, s \in S, \\
    & \sum_{m \in M} o_{j,k,m,t,s} = p_{j,k,t,s}, \quad \forall j \in J, k \in K, t \in T, s \in S, \\
    & \sum_{m \in M} \lambda_{j,k,m,t,s} = 1, \quad \forall j \in J, k \in K, t \in T, s \in S, \\
    & P_{j,k,t',s} = \sum_{\tau=1}^{7} p_{j,k,(7-\tau)(t'-1)+\tau t',s}, \quad \forall j \in J, k \in K, t' \in T', s \in S,\\
    & \sum_{j \in J} P_{j,k,t',s} + \Slack_{k,t',s} = \Demand_{k,t',s}, \quad \forall k \in K, t' \in T', s \in S, \\
    & \text{Loss}_{t',s} = \sum_{k \in K} \Slack_{k,t',s}, \quad \forall t' \in T', s \in S\\
    & 0 \leq \text{Loss}_{t',s} \leq e^{\max} \sum_{s \in S} (\rho_s \sum_{k\in K} \Demand_{k,t',s}) \delta_{t'}, \quad \forall t' \in T', s \in S, \\
    & \sum_{t' \in T'} \delta_{t'} \leq n, \\
    & \variance = \frac{1}{|T'|} \sum_{t' \in T'} (\mathbb{E}[\text{Loss}_{t',s}^2] - \mathbb{E}[\text{Loss}_{t',s}]^2) = \frac{1}{|T'|} \sum_{t' \in T'} (\sum_{s \in S} \rho_s \text{Loss}_{t',s}^2 - (\sum_{s \in S} \rho_s \text{Loss}_{t',s})^2) \leq l.
\end{align}

Equation (3) indicates that if AI production plant $i$ starts a new batch at time $t$, then the AI product will be produced after $L-1$ time periods (days). In addition, the capacity of each batch is bounded by $[\capacity_i^{\norm},\capacity_i^{\max}]$. Equations (4) and (5) imply that if plant $i$ starts producing the AI at time $t$, the process must be active for the next $L-1$ time periods, during which no new batch can start. Each AI plant $i$ has a cleaning policy to avoid buildup of chemical residues which can adversely impact the purity and yield of subsequent batches. Equations (6) thru (8) ensure that, after producing $B_i$ batches, the plant will undergo cleaning for a total of $C$ periods, during which no new batch may start. The cleaning decision is dictated by the binary variable $\alpha^{\clean}_{i,t,s}$. Equation (9) models the AI plant inventory balance, in which $y_{i,j,t,s}$ denotes the amount of AI product sent from plant $i$ to warehouse $j$ at time $t$ under scenario $s$. Equations (10) thru (12) model the transportation of AI product from plant $i$ to warehouse/distribution center $j$ via mode $m$, in which only one transportation mode can be chosen for each unit of AI product. Note that Equation (10) is nonlinear, and we use McCormick relaxation to obtained a linearized version. Equation 13 indicates that if warehouse $j$ is active at time period $t$, the AI can be transported into and/or out of it. Similar to Equation (9), Equation (14) performs mass balance for each warehouse/distribution center, where $p_{j,k,t,s}$ is the amount of AI product shipped from warehouse $j$ to the market region $k$ at time $t$ under scenario $s$. Equations (15) thru (17) model the transportation of AI product from warehouse $j$ to market region $k$. And Equation (15) is also linearized using McCormick relaxation. Equation (17) converts daily shipment amount of AI product into weekly amount, which is compared with the weekly demand under each scenario $s$ to calculate possible demand loss. The demand loss for each market region at the end of week $t'$ under scenario $s$ is defined in Equation (19) in terms of the slack variable $\Slack$. There is a limitation in the number of time periods that demand may not be satisfied completely. Constraint $(20)$ shows that the total amount of possible unsatisfied demand at week $t'$ under scenario $s$ must be bounded. And $\text{loss}_{t',s}$ denotes the overall demand loss for all market regions in week $t'$. To ensure robustness and resilience of supply chain, we require the variance of $\text{loss}_{t',s}$ across all planning periods (weeks) to be bounded. Equation (23) is a quadratic constraint, making this problem a MIQCP. In the next section, we replace this nonlinear constraint with its perspective cuts to reformulate the MIQCP into a MILP. 

\subsection{Reformulation Using Perspective Cuts}
To replace Equation (23) with its perspective cuts, we need to first calculate the gradient of variance. From Equation (23), we have:

\begin{equation}
    \frac{\partial \variance}{\partial \text{Loss}_{t'',s'}} = \frac{1}{|T'|}(2 \rho_{s'} \text{Loss}_{t'',s'} - 2 \rho_{s'} \sum_{s \in S} \rho_s \sum_{t' \in T'} \text{Loss}_{t',s}), \quad \forall t'' \in T', s' \in S.
\end{equation}

Thus, gradient of variance can be determined as:
\begin{align*}
    \nabla \variance = \frac{1}{|T'|}\left(\frac{\partial \variance}{\partial \text{Loss}_{1,1}}, \cdots, \frac{\partial \variance}{\partial \text{Loss}_{|T'|,|S|}}\right)^T = \frac{1}{|T'|} & (2 \rho_{1} \text{Loss}_{1,1} - 2 \rho_{1} \sum_{s \in S} \rho_s \sum_{t' \in T'} \text{Loss}_{t',s}, \cdots, \\ 
    & 2 \rho_{|S|} \text{Loss}_{|T'|,|S|} - 2 \rho_{|S|} \sum_{s \in S} \rho_s \sum_{t' \in T'} \text{Loss}_{t',s} )^T
\end{align*}

One can verify that the perspective cut of Equation (23) is:

\begin{equation}
    \begin{split}
        \frac{1}{|T'|} &(2\sum_{t' \in T'} (\sum_{s \in S} \rho_s \text{Loss}_{t',s}^{*}\text{Loss}_{t',s}) - 2 \sum_{s' \in S} \rho_s (\sum_{t'' \in T} \text{Loss}_{t'',s'}) (\sum_{s \in S} \rho_s \sum_{t' \in T'} \text{Loss}_{t',s}^{*}) \\
         & + (\sum_{t' \in T'} (\sum_{s \in S} \rho_s \text{Loss}_{t',s}^{*2} - (\sum_{s \in S} \rho_s \text{Loss}_{t',s}^{*})^2) - l - (2\sum_{t' \in T'} (\sum_{s \in S} \rho_s \text{Loss}_{t',s}^{*2})\\
        & -2 \sum_{s' \in S} \rho_s (\sum_{t'' \in T} \text{Loss}_{t'',s'}^{*}) (\sum_{s \in S} \rho_s \sum_{t' \in T'} \text{Loss}_{t',s}^{*})) \delta_{t}) \leq 0,  
    \end{split}
\end{equation}
which is successively introduced to the reformulated model (MILP) as described earlier.

\section{An Illustrative Case Study}
To illustrate the effectiveness of using perspective cuts to linearize variance constraint, we study a case study involving two AI production plants, two warehouses, and two market regions. The model parameters used for the case study are summarized in the Table \ref{table1}. We consider a total of three different scenarios in the weekly demand for both market regions (see Table \ref{table2}). The original MIQCP model, which contains 1292 continuous variables, 1010 binary variables, and 1344 integer variables, has 5564 linear constraints and 2 nonlinear constraints of Equation (23) (one equality constraint and one inequality constraint). We formulate the problem in JuMP and use SCIP v8.0 in GAMS v40.2.0 to solve the model. We set the solving time to be 600 seconds, at which we obtain a gap of 12.29\% and an objective function value of \S15282.0. For the reformulated MILP model, we only introduce one perspective cut for each planning period (week) to replace the two nonlinear variance constraints. At 600 seconds, we obtain a gap of 8.72\% and an objective function value of \$14886.0, which is 2.59\% less than that of the MIQCP model. We emphasize that the reformulated model always yields a feasible solution to the original MIQCP model, suggesting that it provides a better optimal solution. This is due to the fact that the reformulated model is able to identify a solution with no demand loss, whereas the original MIQCP gives 5 units of demand loss in the second market region during the first week under scenarios 1 and 3. This suggests that the reformulated model produces a more efficient supply chain compared to the original model. It is also worth pointing out that, while adding perspective cuts for nonconvex MINLPs may result in less significant computational time reduction benefits, it will reduce the size of branch-and-cut trees and strengthens the root node relaxation \cite{bestuzhev21}. Having said that, in the case of our MIQCP model, introducing perspective cuts yields a consistent improvement of computational performance, as the original formulation is convex \cite{bestuzhev21}.

\begin{table}[htb]
\caption{Model parameters used in the case study, in which m.u. stands for mass unit.}\label{table1}
\vspace{-0.7cm}
\begin{center}
\begin{tabular}{|c|c|}
\hline
 \textbf{Model parameters} & \textbf{Values} \\ \hline
 Probability of scenario 1, 2, and 3 ($\rho_s$) & 0.35, 0.15, 0.5 \\ \hline
 Fixed cost for AI production plant 1 and 2 (\$$\times 1000$) & 100, 120 \\ \hline
 Variable cost for AI production plant 1 and 2 (\$$\times 1000$/m.u.) & 3, 2.5 \\ \hline
 Holding cost for AI production plant 1 and 2 (\$$\times 1000$/m.u.) & 2.5, 2 \\ \hline
 Shipping cost from AI plant to warehouse for mode 1 and 2 (\$$\times 1000$/m.u.) & 150, 175 \\  \hline 
 Fixed cost for warehouse 1 and 2 (\$$\times 1000$) & 75, 70 \\ \hline
 Holding cost for warehouse 1 and 2 (\$$\times 1000$/m.u.) & 5, 7 \\ \hline 
 Demand loss cost (\$$\times 1000$/m.u.) & 50 \\ \hline
 Normal production capacity per batch (m.u.) for AI plant 1 and 2 & 70, 80 \\ \hline
 Maximum production capacity per batch (m.u.) for AI plant 1 and 2 & 90, 100 \\ \hline
 Initial AI product inventory level (m.u.) for AI plant 1 and 2 & 30, 40\\  \hline
 Minimum shipping capacity (m.u.) for mode 1 and 2 & 20, 25 \\ \hline
 Maximum shipping capacity (m.u.) for mode 1 and 2 & 50, 50 \\ \hline 
 Maximum warehouse capacity (m.u) for warehouse 1 and 2 & 200 \\ \hline
 AI production duration (time periods), $L$ & 3 \\ \hline
 Cleaning duration (time periods), $C$ & 1 \\ \hline
 Initial AI product storage (m.u.) in warehouse 1 and 2 & 10, 10 \\ \hline
\end{tabular}
\end{center}
\end{table}

\begin{table}[htb] 
\caption{Scenario-based weekly demand (in m.u.) for each market region used in the case study.}\label{table2}
\vspace{-0.7cm}
\begin{center}
\begin{tabular}{|c|c|c|c|c|}
\hline
 \multirow{2}{*}{\textbf{Scenario}} &  \multicolumn{2}{|c|}{\textbf{Week 1}} & \multicolumn{2}{|c|}{\textbf{Week 2}} \\
 & \textbf{Market 1} & \textbf{Market 2} & \textbf{Market 1} & \textbf{Market 2} \\ \hline
 1 & 100 & 75 & 250 & 125 \\ \hline
 2 & 100 & 75 & 120 & 150 \\ \hline
 3 & 100 & 75 & 110 & 175 \\ 
\hline
\end{tabular}
\end{center}
\end{table}

\section{Conclusion}
In this work, we optimize the supply chain of an agrochemical active ingredient by formulating and solving a scenario-based stochastic MIQCP problem. The nonlinearity of the model comes from the variance constraint that is used to quantify risks associated with unforeseen demand loss. For the first time, we propose to reformulate the variance constraint using perspective reformulation. The reformulated model, which is a MILP, always gives a feasible solution to the original MINLP model. Using a simple case study, we demonstrate the effectiveness of perspective cuts in fostering convergence and reducing computation time. Specifcally, we show that, upon successive introduction of perspective cuts, the optimal solution of the MILP tightly approximate the optimal solution of the original MIQCP.

\end{document}